\documentclass[11pt]{amsart}
\usepackage{amsmath, amssymb, amsthm, mathtools}
\usepackage[margin=1.15in]{geometry}
\usepackage{microtype}
\usepackage{enumitem}
\usepackage{xcolor}
\usepackage[colorlinks=true,linkcolor=blue!50!black,citecolor=blue!50!black,
            urlcolor=blue!50!black]{hyperref}
\usepackage[capitalise,noabbrev]{cleveref}

\newtheorem{theorem}{Theorem}[section]
\newtheorem{lemma}[theorem]{Lemma}
\newtheorem{corollary}[theorem]{Corollary}
\newtheorem{proposition}[theorem]{Proposition}
\theoremstyle{definition}
\newtheorem{definition}[theorem]{Definition}
\newtheorem{example}[theorem]{Example}
\newtheorem{assumption}[theorem]{Standing Assumption}
\theoremstyle{remark}
\newtheorem{remark}[theorem]{Remark}
\crefname{assumption}{Standing Assumption}{Standing Assumptions}
\Crefname{assumption}{Standing Assumption}{Standing Assumptions}
\numberwithin{equation}{section}

\newcommand{\R}{\mathbb{R}}
\newcommand{\K}{\mathcal{K}}
\newcommand{\A}{\mathcal{A}}
\newcommand{\Bor}{\mathcal{B}}
\newcommand{\Sph}{\mathbb{S}^{d-1}}
\newcommand{\Ball}{\mathbb{B}}
\newcommand{\DPhi}{D_{\Phi}}
\newcommand{\Du}[1]{D_{\Phi}^{#1}}
\newcommand{\Dup}{D_{\Phi}^{u}}
\newcommand{\Dupn}[1]{\bigl(D_{\Phi}^{u}\bigr)^{#1}}
\newcommand{\uDup}{\overline{D}_{\Phi}^{\,u}}
\newcommand{\lDup}{\underline{D}_{\Phi}^{\,u}}
\newcommand{\norm}[1]{\left\|#1\right\|}
\newcommand{\abs}[1]{\left|#1\right|}
\newcommand{\inn}[2]{\langle #1,\, #2\rangle}
\newcommand{\dd}{\,\mathrm{d}}
\DeclareMathOperator*{\essinf}{ess\,inf}

\begin{document}

\title[Differentiation with respect to set-valued probabilities]
{Differentiation with respect to set-valued probabilities}

\author{Davide La Torre}
\address{SKEMA Business School, Universit\'e C\^ote d'Azur,
Sophia Antipolis, France}
\email{davide.latorre@skema.edu}
\address{Department of Mathematics, University of Bologna, Bologna, Italy}
\email{davide.latorre3@unibo.it}

\author{Franklin Mendivil}
\address{Department of Mathematics and Statistics, Acadia University,
Wolfville, NS, Canada}
\email{franklin.mendivil@acadiau.ca}
\thanks{The research of F. Mendivil is partially supported by NSERC\ 2019:05237.}

\author{Matteo Rocca}
\address{Department of Economics, University of Insubria, Varese, Italy}
\email{matteo.rocca@uninsubria.it}


\begin{abstract}
We study the differentiation of a real function on $[0,1]$ with respect
to a set-valued probability: a map $\Phi$ assigning to each Borel set a
compact convex body in $\R^{d}$, additive for Minkowski addition and
normalised so that $\Phi([0,1])$ is the closed unit ball. The support
function of $\Phi$ in a direction $u$ is a probability measure
$\Phi_{u}$; differentiating against $\Phi_{u}$ gives a directional
derivative, and averaging over the sphere gives a global derivative.
We show that monotonicity of $\Phi$ follows from the other axioms,
characterise the direction-indexed families of measures that arise,
and identify directional derivates with Dini derivates through a time
change. We prove linearity, product, quotient and chain rules (the last
with a clock-density factor), sharp subadditivity of the upper
derivative, and a scalar reduction of the global derivative with a
harmonic-mean density. We prove one-point and integral mean value
theorems, Cauchy, Darboux and Taylor theorems, and show that a
direction-independent mean point need not exist.
\end{abstract}

\maketitle

\section{Introduction}\label{sec:intro}

Differentiation with respect to a measure replaces the length of an
interval in a difference quotient by its mass under a fixed Borel
measure. In one dimension the \emph{$\mu$-derivative} of $f$ at $x$ is
\[
  D_{\mu}f(x)\;:=\;\lim_{y\to x}\frac{f(y)-f(x)}{\mu((x,y])},
  \qquad \mu((x,y]):=-\mu((y,x])\ \text{for }y<x,
\]
which is the ordinary derivative when $\mu$ is Lebesgue measure. The
idea goes back to Lebesgue's differentiation theorem \cite{Lebesgue1910};
Besicovitch \cite{Besicovitch1945,Besicovitch1946} extended it to
general Radon measures by his covering principle, obtaining the
Radon--Nikodym derivative $\dd\nu/\dd\mu$ as the almost-everywhere limit
of the ratios $\nu(I)/\mu(I)$ over shrinking intervals when
$\nu\ll\mu$. Classical derivates and covering arguments are discussed
by Thomson \cite{Thomson1994}. Measure-geometric differentiation on
fractals was developed by Freiberg and Z\"ahle \cite{Freiberg2002}, and
the reduction of the associated harmonic calculus to classical calculus
through distribution functions by Kesseb\"ohmer, Samuel and Weyer
\cite{KessebohmerEtAl2016}; related, but distinct, constructions
include Laplacians on self-similar sets \cite{Kigami1993} and
differentiable structures on doubling metric measure spaces supporting
a Poincar\'e inequality \cite{Cheeger1999}. A recurring theme is that
the reference measure, not the ambient metric, dictates the analytic
behaviour of the derivative.

This paper treats the case in which the reference measure is itself
uncertain and is encoded by a \emph{set-valued probability}: a map
$\Phi$ assigning to each Borel subset of $[0,1]$ a compact convex body
in $\R^{d}$, additive for Minkowski addition and normalised so that
$\Phi([0,1])$ is the closed unit ball. Integration of set-valued maps
goes back to Aumann \cite{Aumann1965} and Hiai and Umegaki
\cite{HiaiUmegaki1977}; limit theorems for random compact sets to
Artstein and Vitale \cite{ArtVit1975} and Puri and Ralescu
\cite{PuriRal1985}; see also \cite{CastVal1977,Molchanov2005}. In
imprecise probability \cite{Walley1991,AugustinEtAl2014}, a credal set
is a convex set of probability measures whose pointwise lower and
upper envelopes define lower and upper probabilities; related
constructions include the multivalued-mapping model of Dempster
\cite{Dempster1967}, the capacities of Choquet \cite{Choquet1954}, the
interval and $p$-box models of \cite{FersonGinz1996,Weich2000}, and the
sublinear expectations of Peng \cite{Peng2019}. The model used here is
geometric: evaluating the support function of $\Phi$ in a direction
$u$ of the unit sphere produces an \emph{additive} probability measure
$\Phi_{u}$, so a set-valued probability is a coherent,
direction-indexed bundle of classical probability measures. It is
therefore distinct from a general lower/upper envelope, which need not
be additive, and it does not represent every credal set.

Our construction is also orthogonal to the derivatives \emph{of}
set-valued maps, such as the contingent derivative \cite{Aubin1981},
the graphical derivative \cite{RockWets1998} and the Hukuhara
derivative \cite{Hukuhara1967}. We do not differentiate $\Phi$; we use
$\Phi$ as the reference measure in the denominator when differentiating
a scalar function $f$. The device that makes this possible is the
support-function embedding $K\mapsto h_{K}$, an isometry from compact
convex bodies with the Hausdorff metric into $C(\Sph)$
\cite{Rockafellar1970,Schneider2014,Hormander1994}: in each fixed
direction it turns Minkowski arithmetic into scalar arithmetic, so the
apparatus of differentiation with respect to a measure applies
direction by direction.

\subsection*{Main results}
The normalisation $\Phi([0,1])=\Ball^{d}$ makes each directional slice
$\Phi_{u}$ a probability measure (\cref{prop:slice}); a family
$\{\nu_{u}\}$ of probability measures arises from a set-valued
probability exactly when $u\mapsto\nu_{u}(A)$ is a support function for
every Borel set $A$ (\cref{prop:representation}); and set-monotonicity
of $\Phi$ is a consequence of finite additivity and the condition
$0\in\Phi(A)$, not an independent axiom (\cref{prop:monotone}).

The \emph{directional derivative} $\Dup f$ is formed with the signed
increment of the distribution function $F_{u}$ of $\Phi_{u}$, and the
\emph{global derivative} $\DPhi f$ is its average over the sphere. The
time change $t=F_{u}(x)$ identifies the directional upper and lower
derivatives with classical Dini derivates (\cref{prop:dini}), which
gives almost-everywhere differentiability of monotone functions in each
direction. Under positivity and integrability hypotheses the global
derivative reduces to differentiation with respect to a single scalar
measure, with a harmonic-mean density in the Lebesgue case
(\cref{thm:scalar-reduction,cor:harmonic}); the directional family,
rather than its average, retains the separate clocks.

We prove linearity and product and quotient rules, and show that the
upper derivative is only subadditive, with equality when one summand is
differentiable and strict inequality in an explicit example
(\cref{prop:upper-additive,ex:strict}). The chain rule carries a
clock-density factor $\rho_{u}(g(x))/\rho_{u}(x)$, which is necessary
outside the Lebesgue case (\cref{thm:chain,ex:chain}). The mean value
theorem holds in a one-point form in each direction and in an integral
form; the intermediate point depends on the direction and a single
direction-independent point need not exist
(\cref{thm:mvt,rem:no-single-point,ex:two-clocks}). We also obtain
constancy and monotonicity criteria, a Darboux property, a Cauchy mean
value theorem, a Jensen inequality, the differentiation of indefinite
integrals, and a directional Taylor formula with Lagrange remainder
(\cref{thm:taylor}). Worked examples---the interval-valued model, a
genuinely anisotropic planar family, a corner, and the Cantor
measure, which marks the boundary of our regularity assumption---are
placed alongside the results they illustrate.

\subsection*{Organisation}
\Cref{sec:prelim} recalls support functions and gives the axiomatics of
set-valued probabilities. \Cref{sec:def} introduces the directional and
global derivatives, relates them to Dini derivates and proves the
scalar reduction. \Cref{sec:properties} establishes the algebraic rules
and the chain rule, and \Cref{sec:mvt} the mean value theorems and
their consequences. \Cref{sec:conclusion} contains concluding remarks.

\section{Preliminaries}\label{sec:prelim}

\subsection{Notation}
Throughout, $\Omega=[0,1]$, $\A=\Bor([0,1])$ is its Borel
$\sigma$-algebra and $\lambda$ is Lebesgue measure on $\Omega$. We fix a
dimension $d\ge 1$ and write $\K(\R^{d})$ for the family of nonempty
compact convex subsets of $\R^{d}$, with the Hausdorff metric
\[
d_{H}(K,L)=\max\Bigl\{\sup_{x\in K}\inf_{y\in L}\norm{x-y},\,
\sup_{y\in L}\inf_{x\in K}\norm{x-y}\Bigr\}.
\]
The closed unit ball is
$\Ball^{d}$ and $\Sph$ is the unit sphere. We fix a Borel probability
measure $\mu$ on $\Sph$; unless stated otherwise $\mu$ is the normalised
surface measure (for $d=1$, $\Sph=\{+1,-1\}$ and
$\mu=\tfrac12(\delta_{+1}+\delta_{-1})$; for $d=2$,
$\mathrm{d}\mu=\mathrm{d}\theta/2\pi$).

\subsection{Support functions}

\begin{definition}[Support function]\label{def:supp}
For $K\in\K(\R^{d})$ and $u\in\R^{d}$, the \emph{support function} of
$K$ is $h_{K}(u):=\sup_{v\in K}\inn{u}{v}$.
\end{definition}

\begin{proposition}[{\cite[Secs.~1.7--1.8]{Schneider2014}}]
\label{prop:supp-props}
For $K,L\in\K(\R^{d})$, $u\in\R^{d}$ and $c\ge 0$:
\begin{enumerate}[label=\textup{(\roman*)}, leftmargin=2.4em]
  \item $h_{K}$ is convex and positively homogeneous, hence sublinear,
        and $h_{cK}=c\,h_{K}$;
  \item $h_{K+L}=h_{K}+h_{L}$, where $K+L$ is the Minkowski sum;
  \item $K\subseteq L$ if and only if $h_{K}(u)\le h_{L}(u)$ for all
        $u\in\Sph$;
  \item $d_{H}(K,L)=\sup_{u\in\Sph}\abs{h_{K}(u)-h_{L}(u)}$, so
        $K\mapsto h_{K}|_{\Sph}$ is an isometry into $C(\Sph)$;
  \item $h_{\Ball^{d}}(u)=\norm{u}$; in particular
        $h_{\Ball^{d}}(u)=1$ for $u\in\Sph$.
\end{enumerate}
\end{proposition}

\subsection{Set-valued probabilities}

\begin{definition}[Set-valued probability]\label{def:svp}
A \emph{set-valued probability} on $(\Omega,\A)$ with values in
$\R^{d}$ is a map $\Phi:\A\to\K(\R^{d})$ such that
\begin{enumerate}[label=\textup{(P\arabic*)}, leftmargin=3em]
  \item\label{P1} \emph{normalisation}: $\Phi(\Omega)=\Ball^{d}$;
  \item\label{P2} \emph{finite additivity}: $\Phi(A\cup B)=\Phi(A)+
        \Phi(B)$ for disjoint $A,B\in\A$ (Minkowski sum);
  \item\label{P3} \emph{origin condition}: $0\in\Phi(A)$ for all
        $A\in\A$;
  \item\label{P4} \emph{continuity from below}: if $A_{n}\uparrow A$
        in $\A$, then $d_{H}\big(\Phi(A_{n}),\Phi(A)\big)\to 0$.
\end{enumerate}
\end{definition}

Set-monotonicity need not be imposed as a separate axiom.

\begin{proposition}[Monotonicity]\label{prop:monotone}
If $\Phi$ satisfies \ref{P2} and \ref{P3}, then $A\subseteq B$ implies
$\Phi(A)\subseteq\Phi(B)$.
\end{proposition}

\begin{proof}
Write $B=A\sqcup(B\setminus A)$. By \ref{P2},
$\Phi(B)=\Phi(A)+\Phi(B\setminus A)$, and since
$0\in\Phi(B\setminus A)$ by \ref{P3},
$\Phi(A)=\Phi(A)+\{0\}\subseteq\Phi(A)+\Phi(B\setminus A)=\Phi(B)$.
\end{proof}

For $u\in\Sph$ define the \emph{directional support measure}
\begin{equation}\label{eq:Phiu}
  \Phi_{u}(A)\;:=\;h_{\Phi(A)}(u)\;=\;\sup_{v\in\Phi(A)}\inn{u}{v},
  \qquad A\in\A.
\end{equation}

\begin{proposition}[Directional slices]\label{prop:slice}
For every $u\in\Sph$, the set function $\Phi_{u}$ in \eqref{eq:Phiu} is
a Borel probability measure on $\Omega$.
\end{proposition}

\begin{proof}
\emph{Finite additivity.} For disjoint $A,B$,
$\Phi_{u}(A\cup B)=h_{\Phi(A)+\Phi(B)}(u)=h_{\Phi(A)}(u)+h_{\Phi(B)}(u)$
by \ref{P2} and \cref{prop:supp-props}(ii).

\emph{Empty set.} Applying \ref{P2} to $\Omega=\Omega\sqcup\emptyset$
gives $h_{\Phi(\Omega)}=h_{\Phi(\Omega)}+h_{\Phi(\emptyset)}$, hence
$h_{\Phi(\emptyset)}\equiv0$ and $\Phi_{u}(\emptyset)=0$ (indeed
$\Phi(\emptyset)=\{0\}$ by \ref{P3} and \cref{prop:supp-props}(iii)).

\emph{Nonnegativity.} By \ref{P3},
$\Phi_{u}(A)=\sup_{v\in\Phi(A)}\inn{u}{v}\ge\inn{u}{0}=0$; monotonicity
is not needed.

\emph{Normalisation.} $\Phi_{u}(\Omega)=h_{\Ball^{d}}(u)=1$ by \ref{P1}
and \cref{prop:supp-props}(v).

\emph{Countable additivity.} Let $\{B_{k}\}$ be disjoint with union
$B$ and put $A_{n}=\bigcup_{k\le n}B_{k}\uparrow B$. Then
$\Phi_{u}(A_{n})=\sum_{k\le n}\Phi_{u}(B_{k})$, while \ref{P4} and
\cref{prop:supp-props}(iv) give
$\abs{\Phi_{u}(A_{n})-\Phi_{u}(B)}\le d_{H}(\Phi(A_{n}),\Phi(B))\to0$.
Hence $\Phi_{u}(B)=\sum_{k}\Phi_{u}(B_{k})$.
\end{proof}

Thus a set-valued probability is a \emph{bundle}
$\{\Phi_{u}\}_{u\in\Sph}$ of Borel probability measures on $[0,1]$,
glued together by the support-function representation. The next result
characterises the bundles that arise and shows that the axioms are
consistent.

\begin{proposition}[Representation of bundles]\label{prop:representation}
Let $\{\nu_{u}\}_{u\in\Sph}$ be Borel probability measures on
$\Omega$. There is a set-valued probability $\Phi$ with $\Phi_{u}=\nu_{u}$
for all $u$ if and only if
\begin{enumerate}[label=\textup{(\alph*)}, leftmargin=2.2em]
\item for every $A\in\A$, the map $u\mapsto\nu_{u}(A)$ is the
  restriction to $\Sph$ of a support function, i.e.\ its positively
  $1$-homogeneous extension $z\mapsto\norm{z}\,\nu_{z/\norm{z}}(A)$,
  defined to be zero at $z=0$, is sublinear on $\R^{d}$; and
\item $\nu_{u}(\Omega)=1$ for all $u$.
\end{enumerate}
In that case $\Phi$ is unique: $\Phi(A)$ is the body whose support
function is given in \textup{(a)}.
\end{proposition}

\begin{proof}
\emph{Necessity.} If $\Phi$ exists, $\nu_{u}(A)=h_{\Phi(A)}(u)$ gives
(a), and $\nu_{u}(\Omega)=h_{\Ball^{d}}(u)=1$ gives (b).

\emph{Sufficiency.} By the characterisation of support functions as the
sublinear functions \cite[Thm.~1.7.1]{Schneider2014}, (a) determines a
unique $\Phi(A)\in\K(\R^{d})$ with $h_{\Phi(A)}(u)=\nu_{u}(A)$. For
\ref{P1}, $h_{\Phi(\Omega)}=1=h_{\Ball^{d}}$ on $\Sph$. For \ref{P2},
additivity of each $\nu_{u}$ gives $h_{\Phi(A\cup B)}=h_{\Phi(A)}+
h_{\Phi(B)}=h_{\Phi(A)+\Phi(B)}$, whence $\Phi(A\cup B)=\Phi(A)+\Phi(B)$
by \cref{prop:supp-props}(iii). For \ref{P3}, $h_{\Phi(A)}\ge0$, so
$0\in\Phi(A)$. For \ref{P4}, if $A_{n}\uparrow A$ then
$\nu_{u}(A_{n})\uparrow\nu_{u}(A)$ for each $u$; these are continuous
functions of $u$ (they are support functions) increasing to a continuous limit
on the compact set $\Sph$, so by Dini's theorem the convergence is
uniform, which by \cref{prop:supp-props}(iv) is
$d_{H}(\Phi(A_{n}),\Phi(A))\to0$.
\end{proof}

\begin{example}[Isotropic and interval families]\label{ex:nonvacuous}
For a single probability measure $\nu$, the constant bundle
$\nu_{u}\equiv\nu$ satisfies (a)--(b), since $z\mapsto\nu(A)\norm{z}$ is
sublinear; this gives the \emph{isotropic} set-valued probability
$\Phi(A)=\nu(A)\Ball^{d}$ (\cref{ex:lebesgue,ex:cantor}). For $d=1$,
every pair $\nu_{+1}=P_{+}$, $\nu_{-1}=P_{-}$ satisfies (a)--(b), since
every nonnegative function on $\mathbb S^{0}=\{\pm1\}$ is a support
function; this is the interval model of \cref{ex:interval}.
\end{example}

For each $u$ we write the distribution function
\begin{equation}\label{eq:cdf}
  F_{u}(x)\;:=\;\Phi_{u}\big([0,x]\big),\qquad x\in[0,1],
\end{equation}
which is nondecreasing and right-continuous with $F_{u}(1)=1$. For
$0\le x\le y\le1$, $F_{u}(y)-F_{u}(x)=\Phi_{u}((x,y])$; we use the signed
increment $F_{u}(y)-F_{u}(x)$ for $x,y$ in either order.

It is also interesting to notice that there is a similar representation of set-valued probabilities in terms of  directional density functions.
As the proof of this result is similar to that of Proposition \ref{prop:representation} we omit the details.

\begin{proposition}
\label{prop:densityrepresentation}
Let $f_u(t) \ge 0$ be functions continuous for $(u,t) \in \Sph \times [0,1]$.  Then there is a set-valued probability $\Phi$ with $\Phi_u(A) = \int_A f_u(t) \ dt$ for all $u$ if and only if
\begin{enumerate}[label=\textup{(\alph*)}, leftmargin=2.2em]
    \item $\int_0^1 f_u(t) \ dt = 1$; and
    \item for all $t$, the map $u \mapsto f_u(t)$ is the restriction to $\Sph$ of a support function.
\end{enumerate}
In this case, $\Phi(A) = \int_A \phi(t) \ dt$, with $\phi(t) \in \K(\R^d)$ determined by the family $f_u(t)$.
\end{proposition}

\begin{assumption}[Nondegeneracy and regularity]\label{ass:reg}
For $\mu$-almost every $u\in\Sph$, $\Phi_{u}$ is nonatomic and $F_{u}$
is continuous and strictly increasing on $[0,1]$. Such directions are
called \emph{admissible}; we fix once and for all a Borel set of
admissible directions of full $\mu$-measure, and directional statements
are understood on this set unless regularity is explicitly assumed for
every direction. Where a density is needed we assume in addition that
$\Phi_{u}\ll\lambda$, write $\rho_{u}=\dd\Phi_{u}/\dd\lambda$, and assume
$\rho_{u}>0$ $\lambda$-a.e. Whenever densities are integrated over
$(x,u)$ we use jointly measurable versions, which exist for this
dominated Borel kernel; pointwise density hypotheses refer to the
stated versions or to the actual derivative $F_{u}'(x)$.
\end{assumption}

\Cref{ass:reg} holds in the Lebesgue case (\cref{ex:lebesgue}) and for
the anisotropic planar family of \cref{ex:aniso}, but fails for the
Cantor measure (\cref{ex:cantor}). It guarantees that every difference
quotient below has a nonzero denominator and that $t=F_{u}(x)$ is a
change of variables.

\begin{remark}[Time change]\label{rem:cov}
Under \cref{ass:reg}, $F_{u}(0)=0$ and $x\mapsto F_{u}(x)$ is an
increasing homeomorphism of $[0,1]$ onto itself, and $\Phi_{u}$ is the
image of Lebesgue measure under $F_{u}^{-1}$. Hence for every
$g\in L^{1}(\Phi_{u})$ and $[a,b]\subseteq[0,1]$,
$\int_{[a,b]}g\dd\Phi_{u}=\int_{F_{u}(a)}^{F_{u}(b)}(g\circ F_{u}^{-1})(s)
\dd s$. This change of variables reduces $\Phi_{u}$-calculus on $[0,1]$
to Lebesgue calculus on $[0,1]$.
\end{remark}

\begin{remark}
Suppose that $F_u' = \rho_u$ for all $u \in \Sph$.
Then if $\lambda(A) = 0$, we have $\Phi_u(A) = 0$ for all $u$ and so $\Phi(A) = \{ 0 \}$ and thus $\Phi \ll \lambda$.
This means that the set-valued measure $\Phi$ has a set-valued density $\phi$ and $\phi_u = \rho_u$.
Clearly the converse is also true, if $\Phi$ has a set-valued density then $F_u' = \rho_u$ exists for each $u$.
Proposition \ref{prop:densityrepresentation} is the formal expression of this relationship.
\end{remark}

\section{The \texorpdfstring{$\Phi$}{Phi}-derivative}\label{sec:def}

\subsection{Directional and global derivatives}

\begin{definition}[Directional $\Phi$-derivative]\label{def:dir-deriv}
Let $\Phi$ satisfy \cref{ass:reg}, let $u$ be admissible and
$f:[0,1]\to\R$. The \emph{upper} and \emph{lower directional
$\Phi$-derivatives} of $f$ at $x\in(0,1)$ are
\[
  \uDup f(x):=\limsup_{y\to x}\frac{f(y)-f(x)}{F_{u}(y)-F_{u}(x)},
  \qquad
  \lDup f(x):=\liminf_{y\to x}\frac{f(y)-f(x)}{F_{u}(y)-F_{u}(x)}.
\]
When they coincide and are finite, the common value is denoted
$\Dup f(x)$ and $f$ is called \emph{$\Phi$-differentiable at $x$ in the
direction $u$}.
\end{definition}

The denominator is the \emph{signed} increment, which has the sign of
$y-x$. This is essential: with the unsigned mass $\Phi_{u}([x,y])$ the
left-hand difference quotients of an increasing function would be
negative, and the monotonicity results of \cref{sec:mvt} would fail.
Since $\Phi_{u}$ is a probability measure (\cref{prop:slice}), $\Dup f$
is a pointwise derivative with respect to a scalar measure. On the
class of functions admitting the integral representation of
Freiberg and Z\"ahle \cite{Freiberg2002}, it agrees $\Phi_{u}$-a.e.\
with their measure-geometric derivative, by \cref{prop:indefinite}.

\begin{definition}[Global $\Phi$-derivative]\label{def:global-deriv}
The function $f$ is \emph{globally $\Phi$-differentiable at $x$} if it
is $\Phi$-differentiable at $x$ in the direction $u$ for $\mu$-a.e.\ $u$
and $u\mapsto\Dup f(x)$ is $\mu$-integrable; its \emph{global
$\Phi$-derivative} is then
\[
  \DPhi f(x)\;:=\;\int_{\Sph}\Dup f(x)\dd\mu(u).
\]
The upper global derivative $\overline{D}_{\Phi}f(x)$ is defined with
$\uDup$ in place of $\Dup$ whenever the integral exists.
\end{definition}

\begin{proposition}[Measurability]\label{prop:meas}
For each $A\in\A$ the map $u\mapsto\Phi_{u}(A)$ is continuous on
$\Sph$. If $f$ is continuous on a neighbourhood of $x$, then
$u\mapsto\uDup f(x)$ and $u\mapsto\lDup f(x)$ are measurable functions
with values in $\overline{\R}$ on the set of admissible directions.
\end{proposition}

\begin{proof}
A support function of a compact set is convex and finite on $\R^{d}$, hence continuous
\cite[Thm.~10.1]{Rockafellar1970}. Fix $x$ and, for $y\ne x$, put
$q_{y}(u)=(f(y)-f(x))/(F_{u}(y)-F_{u}(x))$. On the set of admissible
directions the denominator is nonzero and continuous in $u$, so $q_{y}$
is continuous there. For fixed admissible $u$, $y\mapsto q_{y}(u)$ is
continuous on a punctured neighbourhood of $x$, so its $\limsup$ and
$\liminf$ as $y\to x$ may be taken along rational $y$. Countable upper
and lower limits of continuous functions are measurable.
\end{proof}

\begin{example}[Interval-valued probability, $d=1$]\label{ex:interval}
Let $P_{+},P_{-}$ be Borel probability measures on $[0,1]$ and
$\Phi(A):=[-P_{-}(A),P_{+}(A)]\subset\R$. Then \ref{P1}--\ref{P3} hold,
$\Phi_{+1}=P_{+}$ and $\Phi_{-1}=P_{-}$; the lower endpoint enters through
$u=-1$, where the support function selects $-\inf$. With
$\mu=\tfrac12(\delta_{+1}+\delta_{-1})$,
\[
  \DPhi f(x)=\tfrac12\bigl(D_{P_{+}}f(x)+D_{P_{-}}f(x)\bigr).
\]
Neither measure is assumed to dominate the other, and
$[-P_{-}(A),P_{+}(A)]$ is a convex-body encoding, not an interval of
probability values.
\end{example}

\subsection{Relation to Dini derivates}

The time change $t=F_{u}(x)$ identifies directional $\Phi$-calculus with
ordinary calculus of a reparametrised function, so that the classical
theory of Dini derivates applies.

\begin{proposition}[$\Phi$-derivates are Dini derivates]\label{prop:dini}
Fix an admissible $u$, put $g:=f\circ F_{u}^{-1}$ and $s:=F_{u}(x)$. Then
\[
  \uDup f(x)=\max\{D^{+}g(s),\,D^{-}g(s)\},\qquad
  \lDup f(x)=\min\{D_{+}g(s),\,D_{-}g(s)\},
\]
where $D^{\pm}$ and $D_{\pm}$ are the right and left upper and lower
Dini derivates of $g$. In particular $f$ is $\Phi$-differentiable at $x$
in the direction $u$ if and only if $g$ is differentiable at $s$, and
then $\Dup f(x)=g'(s)$.
\end{proposition}

\begin{proof}
Since $F_{u}$ is an increasing homeomorphism, $y\to x^{\pm}$ if and only
if $t=F_{u}(y)\to s^{\pm}$, and for $y\ne x$
\[
  \frac{f(y)-f(x)}{F_{u}(y)-F_{u}(x)}=\frac{g(t)-g(s)}{t-s}.
\]
A two-sided upper (lower) limit is the maximum (minimum) of the
one-sided ones, and differentiability means that all four Dini derivates
coincide and are finite.
\end{proof}

\begin{lemma}[Classically differentiable functions]\label{lem:classical}
Fix an admissible $u$ and suppose that $F_{u}$ is differentiable at $x$
with $\rho_{u}(x):=F_{u}'(x)>0$. If $f$ is differentiable at $x$, then
$f$ is $\Phi$-differentiable at $x$ in the direction $u$ and
$\Dup f(x)=f'(x)/\rho_{u}(x)$. Conversely, if $\Dup f(x)$ exists, then
$f'(x)$ exists and equals $\rho_{u}(x)\,\Dup f(x)$.
\end{lemma}

\begin{proof}
For $y\neq x$,
\[
  \frac{f(y)-f(x)}{F_{u}(y)-F_{u}(x)}=\frac{f(y)-f(x)}{y-x}\cdot
  \frac{y-x}{F_{u}(y)-F_{u}(x)},
\]
and the second factor tends to $1/\rho_{u}(x)$. For the converse,
multiply the $\Phi_{u}$-quotient by $(F_{u}(y)-F_{u}(x))/(y-x)\to
\rho_{u}(x)$.
\end{proof}

The argument uses only the differentiability of the distribution
function at $x$; it applies verbatim to any finite Borel measure on
$[0,1]$ whose distribution function has a positive derivative at $x$.

\begin{corollary}[Monotone functions]\label{cor:leb-monotone}
If $f$ is monotone on $[0,1]$, then for each admissible $u$ the
derivative $\Dup f(x)$ exists for $\Phi_{u}$-a.e.\ $x$ and
$\Dup f\in L^{1}(\Phi_{u})$.
\end{corollary}

\begin{proof}
The function $g=f\circ F_{u}^{-1}$ is monotone, so by Lebesgue's theorem
it is differentiable at Lebesgue-a.e.\ $s\in[0,1]$ with $g'\in L^{1}$.
Since $\Phi_{u}$ is the image of Lebesgue measure under $F_{u}^{-1}$
(\cref{rem:cov}), \cref{prop:dini} gives differentiability of $f$ in the
direction $u$ at $\Phi_{u}$-a.e.\ $x$, and
$\int\abs{\Dup f}\dd\Phi_{u}=\int\abs{g'(s)}\dd s<\infty$.
\end{proof}

\begin{remark}\label{rem:fixed-x}
The exceptional $\Phi_{u}$-null set in \cref{cor:leb-monotone} depends
on $u$. Consequently, at a \emph{fixed} point $x$ the existence of
$\Dup f(x)$ for $\mu$-a.e.\ $u$, required by
\cref{def:global-deriv}, does not follow from \cref{cor:leb-monotone}.
It does hold, by \cref{lem:classical}, whenever $f$ is differentiable at
$x$ and, for $\mu$-a.e.\ $u$, $F_{u}$ is differentiable at $x$ with
positive derivative; global differentiability then also requires the
integrability of $u\mapsto\rho_{u}(x)^{-1}$ when $f'(x)\ne0$.
\end{remark}

\begin{example}[A corner]\label{ex:corner}
Let $f(x)=\abs{x-\tfrac12}$ and let $\Phi$ be isotropic,
$\Phi(A)=\nu(A)\Ball^{d}$, where $\nu$ has a continuous positive density
$\rho$. By \cref{prop:dini}, the right and left directional derivates at
$\tfrac12$ are $\pm1/\rho(\tfrac12)$, so
\[
  \uDup f(\tfrac12)=\frac{1}{\rho(\tfrac12)},\qquad
  \lDup f(\tfrac12)=-\frac{1}{\rho(\tfrac12)}.
\]
The clock rescales the corner but does not remove it; for $\rho\equiv1$
these are the one-sided slopes $\pm1$.
\end{example}

\begin{example}[An anisotropic planar family]\label{ex:aniso}
Let $d=2$, $u=(\cos\theta,\sin\theta)$, $0<\varepsilon\le\tfrac13$, and
\[
  \nu_{u}(A):=\int_{A}\rho_{u}(x)\dd x,\qquad
  \rho_{u}(x):=1+\varepsilon\cos(2\theta)\cos(2\pi x).
\]
Each $\nu_{u}$ is a probability measure, since
$\rho_{u}\ge1-\varepsilon>0$ and $\int_{0}^{1}\cos(2\pi x)\dd x=0$.
Writing $G(A):=\int_{A}\cos(2\pi x)\dd x$, so that
$\abs{G(A)}\le\lambda(A)$, we have
$p_{A}(\theta):=\nu_{u}(A)=\lambda(A)+\varepsilon G(A)\cos2\theta$. A
$C^{2}$ function $p$ on the circle is the restriction of a support
function if and only if $p''+p\ge0$ \cite[Sec.~1.7]{Schneider2014}:
the tangential second derivative of the homogeneous extension at radius
$r>0$ is $(p''+p)/r$, and convexity across the origin follows here from
$p\ge0$. Since
$p_{A}''+p_{A}=\lambda(A)-3\varepsilon G(A)\cos2\theta\ge
(1-3\varepsilon)\lambda(A)\ge0$, \cref{prop:representation} yields a
unique set-valued probability $\Phi$ with $\Phi_{u}=\nu_{u}$. It
satisfies \cref{ass:reg}, its densities are bounded below by
$1-\varepsilon$, and its clocks depend on the direction (compare
$\theta=0$ with $\theta=\pi/4$). For $f$ differentiable at $x$,
\cref{lem:classical} gives $\Dup f(x)=f'(x)/\rho_{u}(x)$, and the
integral $\frac1{2\pi}\int_{0}^{2\pi}\frac{\mathrm d\theta}{1+a\cos2\theta}
=(1-a^{2})^{-1/2}$, $\abs a<1$, gives
\[
  \DPhi f(x)=\frac{f'(x)}{\sqrt{1-\varepsilon^{2}\cos^{2}(2\pi x)}} .
\]
\end{example}

\subsection{Scalar reduction of the global derivative}

\begin{theorem}[Scalar reduction]\label{thm:scalar-reduction}
Suppose that $\mu$ has full support on $\Sph$, and define
$m(A):=\int_{\Sph}\Phi_u(A)\dd\mu(u)$ and $F_m(x):=m([0,x])$. Then $m$
is a nonatomic Borel probability measure with continuous, strictly
increasing distribution function, and $\Phi_u\ll m$ for every $u$. Let
$r_u=\dd\Phi_u/\dd m$ be jointly measurable densities, and assume that
at a point $x\in(0,1)$, for $\mu$-a.e.\ $u$,
\begin{equation}\label{eq:ratio}
 \lim_{y\to x}\frac{F_u(y)-F_u(x)}{F_m(y)-F_m(x)}=r_u(x)\in(0,\infty),
 \qquad\text{and}\qquad
 w_m(x):=\int_{\Sph}r_u(x)^{-1}\dd\mu(u)<\infty.
\end{equation}
Then $f$ is globally $\Phi$-differentiable at $x$ if and only if the
finite derivative $D_m f(x)$ exists, and in that case
\begin{equation}\label{eq:scalar-reduction}
 \DPhi f(x)=w_m(x)\,D_m f(x).
\end{equation}
The ratio limit in \eqref{eq:ratio} exists for $m$-a.e.\ $x$ and
$\mu$-a.e.\ $u$. Positivity and reciprocal integrability are additional
hypotheses; for example, $r_u(x)\ge\varepsilon_0>0$ for
$(m\otimes\mu)$-a.e.\ $(x,u)$ ensures them almost everywhere.
\end{theorem}

\begin{proof}
Measurability of $u\mapsto\Phi_u(A)$ follows from \cref{prop:meas}, and
monotone convergence gives countable additivity of $m$. If $m(N)=0$,
the nonnegative continuous function $u\mapsto\Phi_u(N)$ vanishes
$\mu$-a.e., hence everywhere because $\mu$ has full support; so
$\Phi_u\ll m$. Nonatomicity of $m$ follows by integrating the singleton
masses, and strict increase of $F_m$ by integrating the positive masses
of intervals in admissible directions. Jointly measurable
Radon--Nikodym densities are obtained as limits of ratios over a
refining sequence of finite interval partitions, and Tonelli's theorem
gives $\int r_u\dd\mu=1$ $m$-a.e.

By \cref{rem:cov} applied to $m$, the function $F_u\circ F_m^{-1}$ is
the indefinite integral of $r_u\circ F_m^{-1}$, hence absolutely
continuous with derivative $r_u\circ F_m^{-1}$ Lebesgue-a.e.;
Fubini's theorem gives the simultaneous ratio limit for $m$-a.e.\ $x$
and $\mu$-a.e.\ $u$.

At a point satisfying \eqref{eq:ratio}, factor
\[
 \frac{f(y)-f(x)}{F_u(y)-F_u(x)}
 =\frac{f(y)-f(x)}{F_m(y)-F_m(x)}\cdot
  \frac{F_m(y)-F_m(x)}{F_u(y)-F_u(x)}.
\]
If $D_mf(x)$ exists, the directional limits equal $D_mf(x)/r_u(x)$,
which is integrable with integral $w_m(x)D_mf(x)$. Conversely, any
single direction in which $\Dup f(x)$ exists and \eqref{eq:ratio} holds
gives $D_mf(x)=r_u(x)\Dup f(x)$.
\end{proof}

When the directional clocks have Lebesgue densities, the reduction
takes a particularly simple form.

\begin{corollary}[Harmonic-mean density]\label{cor:harmonic}
Let $x\in(0,1)$ and suppose that for $\mu$-a.e.\ $u$ the function
$F_u$ is differentiable at $x$ with $\rho_u(x):=F_u'(x)>0$, and that
\[
 0<w(x):=\int_{\Sph}\rho_u(x)^{-1}\dd\mu(u)<\infty .
\]
If $f$ is differentiable at $x$, then $f$ is globally
$\Phi$-differentiable at $x$ and $\DPhi f(x)=w(x)\,f'(x)$. If these
conditions hold for Lebesgue-a.e.\ $x$, put $\rho_H:=1/w$, the harmonic
mean of the directional densities. Then $\rho_H\le\bar\rho:=
\int_{\Sph}\rho_u\dd\mu$ a.e., $\nu_H(\dd x):=\rho_H(x)\dd x$ is a Borel
measure with $\nu_H([0,1])\le1$ (not necessarily a probability), and
$\DPhi f(x)=D_{\nu_H}f(x)$ at every such $x$ at which the distribution
function of $\nu_H$ has derivative $\rho_H(x)$ and $f$ is
differentiable.
\end{corollary}

\begin{proof}
By \cref{lem:classical}, $\Dup f(x)=f'(x)/\rho_u(x)$ for $\mu$-a.e.\
$u$, which is integrable because $w(x)<\infty$; integration gives
$\DPhi f(x)=w(x)f'(x)$. The inequality between harmonic and arithmetic
means gives $\rho_H\le\bar\rho$, and by Tonelli's theorem and Lebesgue's
theorem on derivatives of monotone functions,
$\int_0^1\bar\rho\dd x=\int_{\Sph}\int_0^1F_u'\dd x\dd\mu\le1$. Finally,
the remark after \cref{lem:classical} gives
$D_{\nu_H}f(x)=f'(x)/\rho_H(x)=w(x)f'(x)$.
\end{proof}

\begin{remark}\label{rem:scalar}
More generally, if the hypotheses of \cref{thm:scalar-reduction} hold
for $m$-a.e.\ $x$, the same conclusion holds with
$\nu_H(\dd x):=w_m(x)^{-1}\dd m(x)$: since $w_m^{-1}\le\int r_u\dd\mu=1$,
the density is bounded, and the ratio
$\nu_H((x,y])/m((x,y])$ tends to $w_m(x)^{-1}$ for $m$-a.e.\ $x$ by
Lebesgue's differentiation theorem in the clock $s=F_m(x)$. The
reduction is therefore an almost-everywhere \emph{scalar} one; it does
not reproduce the directional family by a single clock, and the
directional derivatives retain the separate clocks
(\cref{rem:no-single-point}). In particular, it does not encode all the directional information and the full theory remains multidimensional.
Full support of $\mu$ is used only to
obtain domination in every direction; for a general $\mu$ the same
argument gives domination on its support, which suffices for the
averaged derivative.
\end{remark}

\section{Properties of the \texorpdfstring{$\Phi$}{Phi}-derivative}
\label{sec:properties}

Under the hypotheses of \cref{thm:scalar-reduction}, the global rules
below follow from scalar differentiation. We give direct proofs, which
also cover pointwise situations outside the positivity and
reciprocal-integrability hypotheses of that theorem.

\subsection{Linearity and the upper derivative}

\begin{theorem}[Linearity]\label{thm:linearity}
Suppose that $f$ and $g$ are $\Phi$-differentiable at $x$ in the
direction $u$ for $\mu$-a.e.\ $u$, with $u\mapsto\Dup f(x)$ and
$u\mapsto\Dup g(x)$ integrable. Then for all $\alpha,\beta\in\R$ the
function $\alpha f+\beta g$ is globally $\Phi$-differentiable at $x$
and $\DPhi(\alpha f+\beta g)(x)=\alpha\,\DPhi f(x)+\beta\,\DPhi g(x)$.
\end{theorem}

\begin{proof}
For each such $u$ the difference quotient is linear in the function,
so $\Dup(\alpha f+\beta g)(x)=\alpha\Dup f(x)+\beta\Dup g(x)$; integrate
over $\Sph$.
\end{proof}

Differentiability cannot be dropped from \cref{thm:linearity}: the
upper derivative is in general only subadditive.

\begin{lemma}\label{lem:limsup}
Let $a,b$ be real functions on a punctured neighbourhood of $x$. Then
$\limsup_{y\to x}(a+b)\le\limsup_{y\to x}a+\limsup_{y\to x}b$
whenever the right-hand side is not of the form $\infty-\infty$, with
equality if $\lim_{y\to x}b(y)=L$ exists and is finite.
\end{lemma}

\begin{proof}
Subadditivity is standard. If $b\to L\in\R$, subadditivity applied to
$a=(a+b)+(-b)$ gives $\limsup a\le\limsup(a+b)-L$, and together with
$\limsup(a+b)\le\limsup a+L$ this yields equality in $\overline{\R}$.
\end{proof}

\begin{proposition}[Subadditivity of the upper derivative]
\label{prop:upper-additive}
For all $f,g$ and admissible $u$ such that the right-hand side is not of
the form $\infty-\infty$,
\[
  \uDup(f+g)(x)\;\le\;\uDup f(x)+\uDup g(x).
\]
If $g$ is $\Phi$-differentiable at $x$ in the direction $u$, equality
holds and $\uDup(f+g)(x)=\uDup f(x)+\Dup g(x)$. Consequently, if $g$ is
globally $\Phi$-differentiable at $x$, then
$\overline{D}_{\Phi}(f+g)(x)=\overline{D}_{\Phi}f(x)+\DPhi g(x)$,
provided the upper global derivative of $f$ is defined.
\end{proposition}

\begin{proof}
Apply \cref{lem:limsup} to the difference quotients of $f$ and $g$; for
the global statement, integrate the directional equality over $\Sph$.
\end{proof}

\begin{example}[Strict inequality]\label{ex:strict}
Take $d=1$, $\Phi_{+1}=\Phi_{-1}=\lambda$, $x_0=\tfrac12$,
$f(y)=|y-x_0|\sin(1/(y-x_0))$ for $y\ne x_0$, $f(x_0)=0$, and $g=-f$.
The difference quotients of $f$ and $g$ are
$\pm\operatorname{sgn}(y-x_0)\sin(1/(y-x_0))$, so
$\uDup f(x_0)=\uDup g(x_0)=1$, whereas $\uDup(f+g)(x_0)=0<2$. Both
functions must fail to be differentiable for strict inequality, but
this is not sufficient: equality holds for $g=f$.
\end{example}

\subsection{Product and quotient rules}

\begin{theorem}[Product rule]\label{thm:product}
Under the hypotheses of \cref{thm:linearity}, $fg$ is globally
$\Phi$-differentiable at $x$ and
$\DPhi(fg)(x)=f(x)\,\DPhi g(x)+g(x)\,\DPhi f(x)$.
\end{theorem}

\begin{proof}
Fix an admissible $u$ for which both directional derivatives exist.
Since $\Dup f(x)$ is finite and $F_{u}$ is continuous,
\[
  f(y)-f(x)=\frac{f(y)-f(x)}{F_{u}(y)-F_{u}(x)}\,
  \bigl(F_{u}(y)-F_{u}(x)\bigr)\longrightarrow 0,
\]
so $f$ is continuous at $x$. With $\Delta F=F_{u}(y)-F_{u}(x)\ne0$,
\[
  \frac{f(y)g(y)-f(x)g(x)}{\Delta F}
  =f(y)\,\frac{g(y)-g(x)}{\Delta F}+g(x)\,\frac{f(y)-f(x)}{\Delta F},
\]
which tends to $f(x)\Dup g(x)+g(x)\Dup f(x)$. This is a linear
combination of integrable functions of $u$; integrate over $\Sph$.
\end{proof}

\begin{corollary}[Quotient rule]\label{cor:quotient}
If, in addition, $g(x)\ne0$, then $f/g$ is globally
$\Phi$-differentiable at $x$ and
\[
  \DPhi\!\Bigl(\frac{f}{g}\Bigr)(x)
  =\frac{g(x)\,\DPhi f(x)-f(x)\,\DPhi g(x)}{g(x)^{2}} .
\]
\end{corollary}

\begin{proof}
By continuity of $g$ at $x$ (proof of \cref{thm:product}), $g(y)\ne0$
near $x$, and with $\Delta F$ as above
\[
  \frac{1}{\Delta F}\Bigl(\frac{f(y)}{g(y)}-\frac{f(x)}{g(x)}\Bigr)
  =\frac{1}{g(y)g(x)}\Bigl(g(x)\frac{f(y)-f(x)}{\Delta F}
   -f(x)\frac{g(y)-g(x)}{\Delta F}\Bigr).
\]
Let $y\to x$ and integrate over $\Sph$.
\end{proof}

\subsection{Reduction to the classical derivative}

\begin{theorem}[Lebesgue case]\label{thm:reduction}
If $\Phi(A)=\lambda(A)\,\Ball^{d}$ for all $A\in\A$, then for every $x$
at which $f'(x)$ exists, $\Dup f(x)=f'(x)$ for all $u\in\Sph$, and
$\DPhi f(x)=f'(x)$.
\end{theorem}

\begin{proof}
By \cref{prop:supp-props}(i),(v),
$\Phi_{u}(A)=\lambda(A)h_{\Ball^{d}}(u)=\lambda(A)$, so $F_{u}(x)=x$ and
the difference quotient is the ordinary one. Since $\mu(\Sph)=1$,
$\DPhi f(x)=f'(x)$.
\end{proof}

\begin{example}[Lebesgue case]\label{ex:lebesgue}
For $\Phi(A)=\lambda(A)\Ball^{d}$ we have $\Phi_{u}=\lambda$, $F_{u}(x)=x$
and $\rho_{u}\equiv1$ for all $u$, so \cref{ass:reg} holds. The scalar
reduction is trivial ($m=\lambda$, $w\equiv1$), and the mean value
theorems of \cref{sec:mvt} reduce to the classical ones, with a
direction-independent mean point (\cref{rem:no-single-point}).
\end{example}

\begin{example}[The Cantor measure: the boundary of
\cref{ass:reg}]\label{ex:cantor}
Let $\mu_{C}$ be the Cantor--Lebesgue measure and
$\Phi(A)=\mu_{C}(A)\Ball^{d}$. This is a set-valued probability
(\cref{ex:nonvacuous}) with $\Phi_{u}=\mu_{C}$ for all $u$, and $F_{u}$ is
the Cantor function, which is continuous and nondecreasing but constant
on each removed interval. Hence \cref{ass:reg} fails: the full-interval
quotient is undefined on flat pieces and on one side of each gap
endpoint. Restricting to the Cantor set $\mathcal C$ does not produce
the required homeomorphism, since the two endpoints of each removed
interval have the same clock value. One can instead take limits through
points of $\mathcal C$ with nonzero clock increment. For $f=F_u$ every
such quotient equals $1$, including when a gap endpoint is approached
through $\mathcal C$, although the classical derivative of $f$ vanishes
Lebesgue-a.e.; the integral-based derivative of \cite{Freiberg2002}
also assigns derivative $1$ to $F_u$, $\mu_C$-a.e. These restricted and
integral formulations should be distinguished from the full-interval
definition used here. The results below rely on the time change
$t=F_{u}(x)$ and do not apply to this example as stated; extensions
through generalised inverses are not pursued here.
\end{example}

\subsection{Chain rule}

Composition with $g$ changes the local rate of the clock $F_{u}$, which
produces a density factor absent from the classical statement.

\begin{theorem}[Chain rule]\label{thm:chain}
Let $x\in(0,1)$ and let $g:\R\to\R$ be differentiable at $x$ and a
homeomorphism from a neighbourhood of $x$ contained in $(0,1)$ onto a
neighbourhood of $g(x)$ contained in $(0,1)$. Assume that for $\mu$-a.e.\
$u$ the density $\rho_{u}=F_{u}'$ is continuous and positive near $x$
and near $g(x)$, and that $f$ is $\Phi$-differentiable at $g(x)$ in the
direction $u$. Then $f\circ g$ is $\Phi$-differentiable at $x$ in the
direction $u$, with
\[
  \Dup(f\circ g)(x)=\Dup f\big(g(x)\big)\,
  \frac{\rho_{u}\big(g(x)\big)}{\rho_{u}(x)}\;g'(x),
\]
and, whenever the integrand is $\mu$-integrable,
\[
  \DPhi(f\circ g)(x)=g'(x)\int_{\Sph}\Dup f\big(g(x)\big)\,
  \frac{\rho_{u}\big(g(x)\big)}{\rho_{u}(x)}\dd\mu(u).
\]
In particular, if $\rho_{u}(g(x))=\rho_{u}(x)$ for $\mu$-a.e.\ $u$ (for
instance in the Lebesgue case) and $\DPhi f(g(x))$ exists, then
$\DPhi(f\circ g)(x)=\DPhi f(g(x))\,g'(x)$.
\end{theorem}

\begin{proof}
For $y\ne x$ near $x$, local injectivity gives $g(y)\ne g(x)$, and
\cref{ass:reg} gives nonzero increments. Factorise
\[
  \frac{f(g(y))-f(g(x))}{F_{u}(y)-F_{u}(x)}
  =\underbrace{\frac{f(g(y))-f(g(x))}{F_{u}(g(y))-F_{u}(g(x))}}_{(\mathrm{I})}
  \cdot
  \underbrace{\frac{F_{u}(g(y))-F_{u}(g(x))}{g(y)-g(x)}}_{(\mathrm{II})}
  \cdot
  \underbrace{\frac{g(y)-g(x)}{y-x}}_{(\mathrm{III})}
  \cdot
  \underbrace{\frac{y-x}{F_{u}(y)-F_{u}(x)}}_{(\mathrm{IV})}.
\]
As $y\to x$, $g(y)$ ranges over a punctured neighbourhood of $g(x)$, so
$(\mathrm{I})\to\Dup f(g(x))$; moreover $(\mathrm{II})\to\rho_{u}(g(x))$,
$(\mathrm{III})\to g'(x)$ and $(\mathrm{IV})\to1/\rho_{u}(x)$. All four
limits are finite, so no condition on the sign of $g'(x)$ is needed.
Integrate over $\Sph$.
\end{proof}

\begin{remark}[Closed global chain rule]\label{rem:global-chain}
Under the hypotheses of \cref{thm:chain}, suppose that $w(x)$ and
$w(g(x))$ in \cref{cor:harmonic} are finite and positive, and write
$\rho_H=1/w$. Then
\[
 \DPhi(f\circ g)(x)
 =\DPhi f(g(x))\,g'(x)\,\frac{\rho_H(g(x))}{\rho_H(x)}.
\]
Indeed, by \cref{lem:classical} $f$ is differentiable at $g(x)$, and
both sides equal $f'(g(x))g'(x)w(x)$. This is the chain rule for the
scalar reduction.
\end{remark}

\begin{example}[The density factor is necessary]\label{ex:chain}
Let $d=1$ and let $\Phi$ be isotropic with clock $F(x)=x^{2}$, so
$\rho(x)=2x$; \cref{ass:reg} holds although $\rho$ is not bounded away
from zero. Fix $x\in(0,1)$, and let $f(y)=y^{2}$ and $g(x)=\sqrt{x}$, so
$f\circ g(x)=x$ and $\Dup(f\circ g)(x)=1/\rho(x)=1/(2x)$. In the chain
rule, $\Dup f(y)=2y/(2y)=1$, the density ratio is
$\rho(g(x))/\rho(x)=1/\sqrt{x}$ and $g'(x)=1/(2\sqrt{x})$, whose product
is $1/(2x)$. Omitting the density factor would give the incorrect value
$1/(2\sqrt{x})$.
\end{example}

\section{Mean value theorems}\label{sec:mvt}

Throughout this section $f:[a,b]\to\R$ is continuous,
$[a,b]\subseteq[0,1]$ and \cref{ass:reg} holds. For admissible $u$ we
put $\tilde f_{u}:=f\circ F_{u}^{-1}$ on $[F_{u}(a),F_{u}(b)]$, so that
$\tilde f_{u}'(F_{u}(x))=\Dup f(x)$ whenever either side exists
(\cref{prop:dini}).

\begin{definition}\label{def:PhiAC}
The function $f$ is \emph{$\Phi_{u}$-absolutely continuous} on $[a,b]$
if $\tilde f_{u}$ is absolutely continuous on $[F_{u}(a),F_{u}(b)]$.
\end{definition}

If $f$ is $\Phi_{u}$-absolutely continuous, then $\tilde f_{u}$ is
differentiable Lebesgue-a.e., so $\Dup f$ exists $\Phi_{u}$-a.e.\ on
$[a,b]$ (as in the proof of \cref{cor:leb-monotone}).

\begin{theorem}[Mean value theorem]\label{thm:mvt}
Let $f:[a,b]\to\R$ be continuous.
\begin{enumerate}[label=\textup{(\alph*)}, leftmargin=2.4em]
\item \emph{(Directional one-point form.)} Suppose that for $\mu$-a.e.\
  $u$ the derivative $\Dup f(x)$ exists and is finite for every
  $x\in(a,b)$. Then for $\mu$-a.e.\ $u$ there is $c_{u}\in(a,b)$ with
  \begin{equation}\label{eq:mvt-point}
    f(b)-f(a)\;=\;\Dup f(c_{u})\,\bigl(F_{u}(b)-F_{u}(a)\bigr)
             \;=\;\Dup f(c_{u})\,\Phi_{u}\bigl((a,b]\bigr).
  \end{equation}
\item \emph{(Integral form.)} Suppose that $f$ is $\Phi_{u}$-absolutely
  continuous on $[a,b]$ for $\mu$-a.e.\ $u$. Then for $\mu$-a.e.\ $u$
  \begin{equation}\label{eq:mvt-int}
    f(b)-f(a)\;=\;\int_{a}^{b}\Dup f(x)\dd\Phi_{u}(x),
  \end{equation}
  and consequently
  $f(b)-f(a)=\int_{\Sph}\int_{a}^{b}\Dup f(x)\dd\Phi_{u}(x)\dd\mu(u)$.
\end{enumerate}
\end{theorem}

\begin{proof}
Fix an admissible $u$ and write $\tilde f=\tilde f_{u}$,
$s_{a}=F_{u}(a)<s_{b}=F_{u}(b)$.

(a) The function $\tilde f$ is continuous on $[s_{a},s_{b}]$ and, by
\cref{prop:dini}, differentiable on $(s_{a},s_{b})$ with
$\tilde f'(t)=\Dup f(F_{u}^{-1}(t))$. The classical mean value theorem
gives $\tau\in(s_{a},s_{b})$ with
$\tilde f(s_{b})-\tilde f(s_{a})=\tilde f'(\tau)(s_{b}-s_{a})$; put
$c_{u}=F_{u}^{-1}(\tau)$.

(b) By the fundamental theorem of calculus for absolutely continuous
functions, $\tilde f(s_{b})-\tilde f(s_{a})=\int_{s_{a}}^{s_{b}}\tilde
f'(t)\dd t$, and the change of variables of \cref{rem:cov} turns the
right-hand side into $\int_{a}^{b}\Dup f\dd\Phi_{u}$. The inner integral
is therefore the constant $f(b)-f(a)$ for $\mu$-a.e.\ $u$, hence
measurable in $u$ with $\mu$-integral $f(b)-f(a)$; no joint
measurability in $(x,u)$ is needed.
\end{proof}

Part (b) has a converse, which completes the directional fundamental
theorem of calculus.

\begin{proposition}[Differentiation of indefinite integrals]
\label{prop:indefinite}
Let $u$ be admissible and $g\in L^{1}(\Phi_{u})$, and put
$G(x):=\int_{[0,x]}g\dd\Phi_{u}$. Then $G$ is continuous and
$\Phi_{u}$-absolutely continuous on $[0,1]$, and $\Dup G(x)=g(x)$ for
$\Phi_{u}$-a.e.\ $x\in(0,1)$.
\end{proposition}

\begin{proof}
By \cref{rem:cov}, $\tilde G(t):=G(F_{u}^{-1}(t))=\int_{0}^{t}\tilde
g(s)\dd s$ with $\tilde g=g\circ F_{u}^{-1}\in L^{1}([0,1])$, so
$\tilde G$ is absolutely continuous and, by Lebesgue's differentiation
theorem \cite{Lebesgue1910}, $\tilde G'(t)=\tilde g(t)$ for
Lebesgue-a.e.\ $t$. By \cref{prop:dini},
$\Dup G(x)=\tilde G'(F_{u}(x))=g(x)$ for $\Phi_{u}$-a.e.\ $x$, since
$F_{u}^{-1}$ carries Lebesgue measure to $\Phi_{u}$.
\end{proof}

\begin{remark}[No direction-independent mean point]
\label{rem:no-single-point}
The point $c_{u}$ in \eqref{eq:mvt-point} depends on the direction,
and in general no single $c\in(a,b)$ satisfies \eqref{eq:mvt-point} for
all $u$ simultaneously; \cref{ex:two-clocks} gives explicit clocks and a
function for which the directional mean points differ. A
direction-independent mean point exists in the isotropic case
$F_{u}\equiv F$ (e.g.\ \cref{ex:lebesgue}), where \eqref{eq:mvt-point}
is a single equation. 
Normally after averaging we only recover the integral form given in \eqref{eq:mvt-int}.
\end{remark}

\begin{corollary}[Lipschitz bound]\label{cor:lipschitz}
If $\abs{\Dup f(x)}\le M$ for every $x\in(a,b)$ and every $u$ in a set
of full $\mu$-measure, then
\[
  \abs{f(b)-f(a)}\;\le\;M\,\essinf_{u\in\Sph}\Phi_{u}\bigl((a,b]\bigr)
  \;\le\;M\int_{\Sph}\Phi_{u}\bigl((a,b]\bigr)\dd\mu(u).
\]
\end{corollary}

\begin{proof}
By \eqref{eq:mvt-point}, $\abs{f(b)-f(a)}\le M\,\Phi_{u}((a,b])$ for
$\mu$-a.e.\ $u$. The left-hand side does not depend on $u$, so we may
take the essential infimum (the infimum if $\mu$ has full support, by
continuity of $u\mapsto\Phi_{u}((a,b])$), and $\essinf\le\int\dd\mu$.
\end{proof}

\begin{corollary}[Constancy]\label{cor:constant}
A continuous $f$ is constant on $[a,b]$ if and only if, for $\mu$-a.e.\
$u$, $\Dup f$ exists and vanishes at every point of $(a,b)$.
\end{corollary}

\begin{proof}
If $f$ is constant, every difference quotient vanishes. Conversely,
apply \eqref{eq:mvt-point} on $[a,x]$ for each $x\in(a,b]$.
\end{proof}

\begin{corollary}[Monotonicity]\label{cor:monotone-f}
Let $u$ be admissible and $f$ be $\Phi_{u}$-absolutely continuous on
$[a,b]$. If $\Dup f(x)\ge0$ at every $x\in(a,b)$ at which it exists,
then $f$ is nondecreasing; if moreover $\Dup f>0$ $\Phi_{u}$-a.e., then
$f$ is strictly increasing.
\end{corollary}

\begin{proof}
The function $\tilde f=f\circ F_{u}^{-1}$ is absolutely continuous with
$\tilde f'(s)=\Dup f(F_{u}^{-1}(s))\ge0$ for Lebesgue-a.e.\ $s$, so it is
nondecreasing, and so is $f=\tilde f\circ F_{u}$. If $\tilde f'>0$
a.e., then $\tilde f$ and $f$ are strictly increasing.
\end{proof}

\begin{proposition}[Darboux property]\label{prop:darboux}
Let $u$ be admissible and suppose that $\Dup f(x)$ exists at every
$x\in(a,b)$. If $a<x_{1}<x_{2}<b$ and $\zeta$ lies strictly between
$\Dup f(x_{1})$ and $\Dup f(x_{2})$, there is $x_{3}\in(x_{1},x_{2})$
with $\Dup f(x_{3})=\zeta$. In particular $\Dup f$ has no jump
discontinuities.
\end{proposition}

\begin{proof}
By \cref{prop:dini}, $\tilde f=f\circ F_{u}^{-1}$ is differentiable on
$(F_{u}(a),F_{u}(b))$ with $\tilde f'(F_{u}(x))=\Dup f(x)$. Darboux's
theorem for $\tilde f'$ gives $s_{3}$ between $F_u(x_1)$ and $F_u(x_2)$
with $\tilde f'(s_{3})=\zeta$; put $x_{3}=F_{u}^{-1}(s_{3})$.
\end{proof}

\begin{theorem}[Cauchy mean value theorem]\label{thm:cauchy}
Let $f,g:[a,b]\to\R$ be continuous and let $u$ be admissible. If
$\Dup f$ and $\Dup g$ exist on $(a,b)$ and $\Dup g$ does not vanish
there, then there is $c_{u}\in(a,b)$ with
\[
  \frac{f(b)-f(a)}{g(b)-g(a)}=\frac{\Dup f(c_{u})}{\Dup g(c_{u})} .
\]
\end{theorem}

\begin{proof}
By \cref{thm:mvt}(a), applied in the direction $u$,
$g(b)-g(a)=\Dup g(c)(F_{u}(b)-F_{u}(a))\ne0$ for some $c$, so the
left-hand side is defined. In the clock $s=F_{u}(x)$ the directional
derivatives of $f$ and $g$ are the classical derivatives of
$\tilde f=f\circ F_{u}^{-1}$ and $\tilde g=g\circ F_{u}^{-1}$, and
$\tilde g'$ does not vanish. The classical Cauchy mean value theorem
gives $\tau\in(F_{u}(a),F_{u}(b))$ with the required identity for
$\tilde f,\tilde g$; put $c_{u}=F_{u}^{-1}(\tau)$.
\end{proof}

\begin{proposition}[Jensen inequality]\label{prop:jensen}
Let $\varphi:\R\to\R$ be convex, $u$ admissible, and $f$
$\Phi_{u}$-absolutely continuous on $[a,b]$. Writing
$\langle h\rangle_{u}:=\Phi_{u}((a,b])^{-1}\int_{a}^{b}h\dd\Phi_{u}$,
\[
  \varphi\bigl(\langle\Dup f\rangle_{u}\bigr)
  \;\le\;\bigl\langle\varphi(\Dup f)\bigr\rangle_{u},
  \qquad
  \langle\Dup f\rangle_{u}=\frac{f(b)-f(a)}{\Phi_{u}((a,b])} .
\]
\end{proposition}

\begin{proof}
The normalised restriction of $\Phi_{u}$ to $(a,b]$ is a probability
measure, so the inequality is Jensen's inequality; the identity is
\eqref{eq:mvt-int} divided by $\Phi_{u}((a,b])$.
\end{proof}

\subsection{A directional Taylor formula}

The mean value theorem is the first-order case of a Taylor expansion in
powers of the increment $F_{u}(y)-F_{u}(x)$. We write $\Dupn{0}f=f$ and
$\Dupn{k}f=\Dup(\Dupn{k-1}f)$ whenever the iterates exist.

\begin{lemma}[Iterated derivatives]\label{lem:iterate}
Let $u$ be admissible and $\tilde f=f\circ F_{u}^{-1}$. If $\tilde f$ is
of class $C^{k}$ on a neighbourhood of $s=F_{u}(x)$, then $\Dupn{j}f$
exists near $x$ for $0\le j\le k$ and $\Dupn{j}f(x)=\tilde
f^{(j)}(F_{u}(x))$.
\end{lemma}

\begin{proof}
For $j=1$ this is \cref{prop:dini}. If $\Dupn{j-1}f=\tilde
f^{(j-1)}\circ F_{u}$ near $x$, then the reparametrisation of
$\Dupn{j-1}f$ is $\tilde f^{(j-1)}$, and \cref{prop:dini} gives
$\Dupn{j}f=\tilde f^{(j)}\circ F_{u}$ near $x$.
\end{proof}

\begin{theorem}[Taylor formula with Lagrange remainder]\label{thm:taylor}
Let $u$ be admissible, $x\ne y$ in $(0,1)$, and suppose that
$\tilde f=f\circ F_{u}^{-1}$ is of class $C^{k+1}$ on an open interval
containing the closed interval between $F_{u}(x)$ and $F_{u}(y)$. Then
there is $\xi$ strictly between $x$ and $y$ with
\[
  f(y)=\sum_{j=0}^{k}\frac{\Dupn{j}f(x)}{j!}\,
       \bigl(F_{u}(y)-F_{u}(x)\bigr)^{j}
   +\frac{\Dupn{k+1}f(\xi)}{(k+1)!}\,
       \bigl(F_{u}(y)-F_{u}(x)\bigr)^{k+1}.
\]
For $k=0$ this is the directional mean value theorem
\eqref{eq:mvt-point}.
\end{theorem}

\begin{proof}
Apply the classical Taylor theorem to $\tilde f$ between $s=F_{u}(x)$
and $\sigma=F_{u}(y)$, obtaining $\tau$ strictly between them, and use
$\tilde f(\sigma)=f(y)$, $\sigma-s=F_{u}(y)-F_{u}(x)$ and
\cref{lem:iterate}, with $\xi=F_{u}^{-1}(\tau)$.
\end{proof}

The second-order term makes $\Dupn{2}f$ a directional
``$\Phi$-Hessian'': if $\tilde f_{u}$ is $C^{2}$ near $F_{u}(x)$,
$\Dup f(x)=0$ and $\Dupn{2}f>0$ near $x$ for a single admissible $u$,
then \cref{thm:taylor} with $k=1$ shows that $x$ is a strict local
minimiser of $f$.

\begin{example}[Direction-dependent mean point]\label{ex:two-clocks}
In the interval model with $d=1$, let $f(x)=x^{2}$ on $[0,1]$, let
$P_{+}$ be Lebesgue measure ($F_{+1}(x)=x$) and let $P_{-}$ have
distribution function $F_{-1}(x)=\sqrt{x}$ (density
$\tfrac12x^{-1/2}$). Then
\[
  \Du{+1}f(x)=2x,
  \qquad
  \Du{-1}f(x)=\frac{2x}{\tfrac12 x^{-1/2}}=4x^{3/2}.
\]
On $[0,1]$, \eqref{eq:mvt-point} reads $1=\Du{u}f(c_{u})$ in both
directions, so
\[
  c_{+1}=\tfrac12,
  \qquad
  c_{-1}=2^{-4/3}\approx0.397 .
\]
No single intermediate point serves both directions, which confirms
\cref{rem:no-single-point}, while the integral form \eqref{eq:mvt-int}
holds in each direction:
$\int_{0}^{1}2x\dd x=1$ and
$\int_{0}^{1}4x^{3/2}\cdot\tfrac12x^{-1/2}\dd x=1$. With
$\mu=\tfrac12(\delta_{+1}+\delta_{-1})$ the global derivative is
$\DPhi f(x)=x+2x^{3/2}$. The iterated derivatives of \cref{lem:iterate}
are $(\Du{+1})^{2}f\equiv2$, since $\tilde f(s)=s^{2}$ in the first
clock, and $(\Du{-1})^{2}f(x)=12x$, since $\tilde f(s)=s^{4}$ in the
second. Both are positive on $(0,1]$, so in each direction $f$ is
convex as a function of the clock variable $s=F_{u}(x)$.
\end{example}

\section{Concluding remarks}\label{sec:conclusion}

The calculus developed here is organised around the bundle
$\{\Phi_{u}\}_{u\in\Sph}$: each slice supplies a one-dimensional
derivative through its distribution function, and the change of
variables $t=F_{u}(x)$ reduces each directional statement to a
classical one. What is new lives in the interaction between different
directions. The global derivative admits a scalar reduction, with a
harmonic-mean density in the Lebesgue case, but the directional family
cannot in general be replaced by a single clock: mean points depend on
the direction, the chain rule carries a clock-density factor, and the
upper derivative is only subadditive. 
The Cantor example indicates the
boundary of the strictly increasing clock hypothesis; extending the
theory to clocks with flat pieces, through generalised inverses or
restricted difference quotients, is a natural next step. 
The stability of
$\DPhi$ under perturbations of $\Phi$ and $f$, and variational problems
built on directional calculus, will be considered elsewhere.

\end{document}